\def\R{I\!\! R}

\def\C{C \!\! \! \!  I\, }
\def\bP{I \! P}
\def\EL{Euler-Lagrange }
\def\Ri{Riemannian }

\def\dd#1#2
{
{{\partial #1}   \over {\partial #2}}
}

\def\ddt#1{
{{d#1}   \over {dt}}}

{\bf INFINITELY MANY ECLIPSES}

\bigskip
\centerline{Richard Montgomery$^{**}$ }
\medskip
\centerline{$^{**}$ Mathematics Dept. UCSC, Santa Cruz, CA 95064 USA}
\smallskip
\centerline{rmont@math.ucsc.edu}
\vskip0.7cm

{\bf Abstract.}{\it
We show  that any
bounded zero-angular  momentum solution
for the Newtonian three-body problem  
must suffer infinitely many eclipses, or collinearities,
provided that it does not suffer a   triple collision.
Motivation for the result comes from the 
dream of building a symbolic dynamics
for the three-body problem, one whose
symbols  1, 2, 3 representing
the three types of eclipses.  The   proof  involves
the conformal geometry of the shape sphere.
}
\vskip .4cm 

1. {\bf Infinitely Many Eclipses.}

A solution to the
 Newtonian three-body problem 
suffers an {\it eclipse} when
 the three bodies, taken to be point masses,  become collinear. 
The solution is  {\it bounded} if
the distances between    bodies
remains   bounded by a fixed constant for all time.

\proclaim Theorem 1.   Every bounded solution  
of the three body problem with   zero angular momentum
and no triple collisions suffers infinitely many eclipses.
 
  Mark Levi conjectured
this theorem during  a conversation with the
author in 1998.   

The Lagrange
solutions show that the theorem fails if we omit
the zero angular momentum condition.   
In these  solutions the  three bodies
form an equilateral  triangle  at every instant.
Bounded  Lagrange solutions with
non-zero angular momentum exist for all time,
and for all mass distributions.
They suffer no eclipses, nor triple
collisions.

The theorem allows  binary collisions
 in which case we use  Levi-Civita  regularization 
to analytically continue the solution through the binary collision,
which counts as an eclipse.  The only obstruction
to infinite time existence for a three-body solution is
 triple collision.  As long as the solution suffers
no triple collision, it can be continued analytically in
(regularized)  time.

\vskip .3cm {\bf 2. Motivation.}

Eclipses come in three types,  labelled
 $1$, $2$, and $3$  
depending on the  mass which 
lies between the other two.   
An {\it eclipse sequence}
is  an infinite sequence  in the letters  $1$, $2$, and $3$.  
We may associate to
each collision-free solution   its eclipse sequence.
If the solution is periodic modulo rotations then
its eclipse sequence is periodic.     The free homotopy type of 
a curve which is periodic modulo rotation, whether a solution or not, is
encoded by its periodic eclipse sequence.
Is every free homotopy realized by a collision-free periodic-modulo-rotation solution?
In other words, does every periodic eclipse sequence
arise as the eclipse sequence of some such solution?
Wu-Yi Hsiang asked me this question in 1996.  It
helped lead to the   rediscovery of
the figure eight solution ( Chenciner and Montgomery [2000]),
a solution with eclipse sequence $123123$.
More generally, we can ask is   every infinite eclipse sequence
realized by a solution?  When we attempt to
realize a given eclipse sequence by
the direct method of the calculus of variations,
the solutions we obtain (if any) are   forced to have zero angular momentum.
See Montgomery [1998]. 
This leads us to ask the following closely related   questions.
Is the set of   
collinear states a kind of   a   slice for the 
zero-angular momentum  three-body dynamics?
If so, does  this slice  lead to a  symbolic dynamics in the  symbols 1, 2, 
and 3?     Theorem 1 is a partial answer to the slice question 
since it asserts that   every zero angular momentum bounded  orbit without triple collision
 must intersect the  alleged collinear ``slice''    an
infinite number of times.

\vskip .3cm {\bf 3. Intuition and Shape Space}

  {\it Shape space}
is the space of oriented congruence classes
of triangles in the plane.  It  is   homeomorphic
to $\R^3$, but is not isometric to it. (See section 11.) 
We will use spherical coordinates $(R, \phi, \theta)$ on shape space.
  $R$ measures the overall size of the triangle, and  is related to the 
triangle's moment of inertia $I$ (formula in next
section) by $R^2 = I$. The variables $(\phi, \theta)$ coordinatize a two-sphere which we call the
{\it shape sphere} and whose points represent oriented similarity
classes of triangles.  
Any  motion of the three bodies
projects to  the motion of a single  point in this shape space.
When  that motion is a zero  angular momentum solution to 
Newton's equation then this  
shape space motion   is defined by a second-order 
differential equation in shape space which itself 
has the form of a Newton's equations,  but now in shape space. 
Under the homeomorphism of shape space with Euclidean 
three-space,
the set of collinear triangles 
is represented by the   $xy$ plane.   The origin of shape
space represents triple collison.   Within
the collinear plane, and issuing forth from the origin, 
 lie  three rays whose points represent the  
binary collision configurations.  
 The zero angular momentum Newton's 
equation written on shape space  says  that 
the three  binary collision rays exert an
attractive   force 
on  the   moving  point.   Since the rays lie in
the collinear plane,  this
force is always directed towards this
 plane.  Levi conjectured,
arguing from mechanical intuition, 
that the
point is obliged to either oscillate up and down
across the collinear  plane or escape to infinity.

\vskip .3cm {\bf 4. An oscillatory area.}

The proof of theorem 1 is based on  a differential
equation for a certain normalized signed area $z$ 
of the triangle formed by the three bodies, and described by theorem 2 below.
The signed area $\Delta$ of
the  triangle whose vertices are ${\bf x_1, x_2, x_3}$
is 
$$\Delta = {  1\over 2} {\bf n} \cdot ({\bf x_2 - x_1})
\times ({\bf x_3 -x_1})$$
 where ${\bf n}$ is the normal to the plane of the triangle.
Define a normalized signed area   by
$$z = {4 \over \sqrt{3}} {\Delta \over {I_1}}$$
where 
$$I_1 = {1 \over 3}(r_{12}^2 + r_{23} ^3 +
r_{31}^2) \hskip .5cm \hbox{with} \hskip .5cm  r_{ij} = |x_i -x_j|
$$
  would  be the moment of inertia of
the triangle  with respect to its
center of mass  {\bf provided}  the masses $m_i$ of its vertices were all $1$.   
$I_1$ is to be compared with the triangle's true moment of
inertia    
 $$I = I_m = (m_1 m_2
r_{12}^2 + m_2 m_3 r_{23} ^3 + m_3 m_1 r_{31}^2)/(m_1 + m_2
+ m_3).$$
  The subscript  $m = (m_1, m_2, m_3)$ indicates the mass
distribution of the three bodies.   There are constants $c, C$
such that $cI_1 \le I \le C I_1$. 
 The motion is bounded if and only if there
is a constant $C_*$ such that 
$I (t) \le C_*$ for all time $t$.  
  The motion
has a triple collision at time $t$ if and only if $I(t) = 0$.
 
 The variable $z$ lies  between $-1$ and $1$,  
with $z = \pm 1$ if and only if the
triangle is Lagrange, i.e. equilateral.  It will be related
to the spherical coordinate $\phi$ mentioned 
briefly in the preceding section by
$z = \sin(\phi)$.  
The solution suffers an eclipse at time $t$ if and only
if $z(t) = 0$.  Thus theorem 1 asserts that $z(t)$ has
infinitely many zeros.

The zero-angular momentum Lagrange solutions,
or {\it Lagrange homothety solutions}  plays a central role in our work here.    
In these solutions an   equilateral
triangle shrinks by homothety to a point in finite time,
thus ending in    
triple collision. 


\proclaim Theorem 2. The normalized
area variable
$z$ satisfies the differential equation
$$\ddt{}(f \dot z) = - q z \hskip 1cm (1) $$
along any zero-angular momentum solution to the three body
problem. The functions 
 $f$ and $q$
are smooth nonnegative  functions, with $f$ a strictly positive 
of shape alone, while     
$q$ is  a function of shape and   velocities which is 
positive  except along initial conditions
for  the  Lagrange homothety solution
where it is zero.  

\noindent Explicit formulae for the    functions $f$ and $q$ of   theorem 2 are
$$f = 3 m_1 m_2 m_3 I_1 ^2/ (m_1 + m_2 + m_3) I  = I \lambda \hskip .5cm (2)$$
and 
$$q  = (1 -  {1 \over 2}{{\cos(\phi)} \over  {\sin(\phi)}}
{1 \over \lambda} {{\partial \lambda} \over {\partial \phi}}
)I (\dot \phi ^2 + \cos^2 (\phi) \dot \theta ^2) 
 - 4{{\cos(\phi)} \over  {\sin(\phi)}} {{\partial U} \over
{\partial \phi}} \hskip .5cm (3),  
$$ 
with
$$\lambda = 3 m_1 m_2 m_3 I_1 ^2/ (m_1 + m_2 + m_3) I^2,
 $$
with $\phi, \theta$  certain spherical   coordinates on the shape
sphere described in section 9,   $\phi$ being related to $z$
by 
$$z = \sin(\phi),$$
and   $U = U(I, \phi, \theta)$ being the negative of the
usual Newtonian potential, viewed as a function on
shape space.

The difficult part of the proof is establishing the
positivity of $q$.  

\proclaim Corollary to the proof of theorem 1 .
The normalized height function $z(t)$ 
of a  zero angular
momentum solution, bounded or not,  has exactly one criticial
point between any two successive zeros, i.e.  successive eclipses,
and this is   a   nondegenerate critical point. In particular,
if the zeros occur at $t_1$ and $t_2$ with $t_1 < t_2$
and if $t_c$ is the critical point, then    $z(t)$
is strictly monotonic on the subintervals 
$t_1 < t <  t_c$
and $t_c < t <  t_2$.

\vskip .3cm {\bf 5. Proof of Theorem 1.}

We prove theorem 1, assuming theorem 2.
An eclipse is a zero of $z$, so we
must show that $z$ has infinitely many zeros.
Equivalently, we show that on any infinite interval
$a \le t \le + \infty$ there is a zero of $z$.

Restrict  attention to the 
case $z(t) >0$.  The argument for $z(t) < 0$ 
proceeds in an  identical manner
except that  the signs of $z$ and its derivative $\dot z$
are   to be reversed. 
 We  first 
show that  if $z(t_1) > 0$ and 
$\dot z (t_1) <0$ then at some later time
$t_2 > t_1$ we must have $z(t_2) = 0$.
Next we will show that if 
$z(t) > 0$  
then eventually for 
some later time
$t_* > t$ we must have $\dot z (t_*) <0$.
Together, these facts  show that $z(t)$
has a zero some finite time later,
and  complete the proof.  

So suppose that  
that   $z(t_1) > 0$ and $\dot z (t_1) < 0$.   
Write $\dot z = {1 \over f} (f \dot z)$
and integrate over the interval $t_1 \le s \le t$
 to obtain 
$$\eqalign{z(t)  &= 
z(t_1) +   \int_{t_1} ^t {1 \over f(s)} (f(s) \dot z
 (s))ds.}$$
Set
$$\delta = - f(t_1) \dot z (t_1),$$
a positive constant.
Since $q \ge 0$ in theorem 2,  differential equation (1),
namely $\ddt{}(f \dot z ) = -q z$, 
says that  
that $f(s) \dot z (s)$ is   monotone decreasing
over any time interval on 
which  $z$ is positive.   That is, $f(s) \dot z (s) < f(t_1)
\dot z (t_1):= - \delta < 0$ for $s > t_1$, as long
as $z(s)$ is positive.  
The boundedness of our
  solution and hence of $I$,  the fact that $\lambda$ is a continuous
positive function on the sphere,
and the fact that $f = I \lambda$ (see eq. (2)) together   imply that
$f$ is bounded.  So there is a positive constant $K$ such that  $0 < f(t) < K$
along our solution. Then $1/f > 1/K$ and
$-1/f < -1/K$.  Consequently 
$\dot z = {(f \dot z)/f}    < - \delta/K$
over our interval of positivity of $z$. 
Now suppose that $z(t)$ remains
positive over the interval $t_1 \le s \le t_2$.  
It follows  
from our integral equation for $z(t)$ and the inequality
immediately above that   
$$
z(t_2)  <   z(t_1) - (\delta/K) (t_2 - t_1).$$
 This inequality together with $z(t) \le 1$ 
forces $z(t_2)$ to be negative 
as soon as
$t_2 -t_1 > K/\delta$.  Consequently
$z$ must have a zero within the time $K/\delta$.

It remains   to show that there
must be a time at which $\dot z$ is negative.
This is equivalent to showing that it is impossible
for a collision-free bounded zero-angular momentum
solution to simultaneously satisfy $z(t) >0$
and $\dot z \ge 0$ over an infinite
time interval $a \le  t < \infty$.  
We argue by contradiction.  Suppose we have
  such a solution.
Since $\dot z \ge  0$ for all $t \ge a$.
the function $z$ is positive
and monotone increasing over
the whole infinite interval, and so
tends to its supremum in infinite positive time.
But $z$ is bounded by $1$, so that we must have   $\dot z
\to 0$.   Again $f = \lambda I$ is bounded. 
It  follows that the limit of $f  \dot z$
as $t \to \infty$ must be zero.
We  now show that the limit
of $\lim_{t \to \infty} z(t) = 1$,
which is to say, that the limiting shape is Lagrange's
equilateral triangle.  
For suppose not. Then $z$ is everywhere
positive and bounded away from Lagrange.
Recall that the coefficient function $q$ 
of the differential equation (1) is non-negative  and
continuous, and is zero if and only if the shape is Lagrange
and the initial conditions are those of Lagrange homothety
solution.  It follows that if   
$\lim_t z(t) < 1$  then 
$q \ge c$ everywhere along our solution,
for some positive constant $c$.
Now use the differential equation (1): 
${d \over dt} (f  \dot
z) = -qz $.  Since $q \ge c > 0$ and
$z > z(a) >0$  the right hand side
of this  differential equation is strictly negative
and bounded away from zero by
the  negative constant $-c z(a)$.  This contradicts  
$\lim_{t \to \infty} f \dot z = 0$.      

Now we know that $z \to 1$ monotonically
as $t \to \infty$ while  
$f \dot z$ decreases monotonically
to zero.    The first fact says the configuration
approaches the Lagrange equilateral shape.  We will now show that
there are times $t_j$ tending to infinity for which the
corresponding 
velocities approach those of the Lagrange homothety
solution.  Integrating the  differential equation (1) 
of theorem 2 from $t = a$
to $\infty$ and using $\lim_{t \to \infty} f(t) \dot z (t) =
0$ 
we obtain
$\int_a ^{\infty} q(s) z(s) ds  = -f(t_1) \dot z (t_1)$. 
  It follows that 
$\int_a ^{\infty} q(s) ds $ 
is finite. This implies that the $\lim \inf$  
of $q$ as $t \to \infty$
is $0$.  Thus  there are time intervals $[t_j, t_{j+1}]$,
$t_j \to \infty$  over which $q(s)$ is
as small as we please.
(We have not excluded the 
possibility   that  
$\lim_{t \to \infty} \sup q(t) > 0.$)
During these intervals of small $q$
the solution is nearly tangent to the Lagrange
homothety configuration, since this is the
only place in phase space where $q$ is zero.
In other words,   the
$\omega$-limit set of our solution curve
contains points of phase space which are 
initial conditions for the   Lagrange homothety
solution.  

It follows that our solution contains
arcs which follow the Lagrange homothety solution 
arbitrarily closely, and hence come
arbitrarily close to the Lagrange 
triple collision.  We now use the results
of Moeckel [1983] on the linearization of 
the flow near Lagrange triple collision.
He performs a McGehee-type blow-up to
add the triple collision states as a boundary
to phase space.  
The  Lagrange triple collision point
becomes a hyperbolic rest point of the resulting
vector field, and the Lagrange homothety solution
lies in its stable manifold.  
We have seen that  our solution
curve comes arbitrarily close to the saddle point,
but does  not lie on its stable manifold, since if it did
 it would suffer a triple
collision.  It follows that the solution curve has
near-collision hyperbolic shaped arcs in which it
closely follows the stable manifold of
the saddle point, coming very close to the point, then makes
a sharp turn  and
follows the unstable manifold to exit
a small neighborhood of the point. 
Consequently its distance in phase space from
the saddle point must decrease.
We will now show that the distance in configuration
space from the Lagrange point must also increase.
Indeed, 
near triple collision the unstable manifold of the
Lagrange point is transverse to the
fibers of the projection $(configuration, velocity) \mapsto
(configuration)$.   This transversality
follows from the same transversality for
the negative  eigenspace of
the linearized flow at Lagrange point.
See Moeckel [1983],
pp. 228-229. Consequently, the  spherical distance 
of our solution from the  Lagrange point must 
 increase. This   distance    can be measured by $1 -z$.
Thus $z$ must decrease hence we must have    
$\dot z <0$ somewhere,  as desired.

QED 

\vskip .3cm{\bf 6. Proof of the Corollary.}
Consider again the case $z > 0$. 
We saw in the proof of theorem 1 that once
$\dot z < 0$ then $z$ continues to decrease monotonically
until it crosses zero. Thus it can have only one
local maximum, on one side of which it is monotone increasing
and the other side of which it is monotone decreasing.
At this maximum we have $\dot z = 0$.
At such a critical point of $z$ 
 eq. (1) of theorem 2 reads
$f \ddot z = - qz$.  It follows that
$\ddot z <0$ at this maximum, since  $f$ 
and $q$ are positive. QED.

\vskip .3cm{\bf 7. Reduced dynamics.}

The proof of theorem 2 boils down to computing
Newton's   equations of
motion for the three bodies using good coordinates  on
 shape space.    Newton's equations 
are the \EL equations for  the Lagrangian  
$$L = {1\over 2} K + U$$
where
$K = m_1 \| \dot x_1 \|^2 + m_2 \|\dot x_2 \|^2
+ m_3 \| \dot x_3 \|^3$
is twice the  kinetic energy,  and 
$U = m_1m_2/r_{12} + m_1 m_2 /r_{13} + m_2 m_3 /
r_{23}$
is the negative of the potential energy.
Here $x_i$, $i =1 ,2, 3$ denote the positions of
the three bodies,  $\dot x_i$ are their velocities,
and   $r_{ij} = \|x_i - x_j \|$ is the distance between
body
$i$ and body $j$.

Shape space is
homeomorphic but not isometric to Euclidean three-space.
Introduce   spherical
coordinates 
$(R,
\phi,
\theta)$ on shape space, with 
$$R^2 = I,$$ 
and $\phi$ being the colatitude,  taken so that
$\phi = 0$ is the equator.  Then  
(Chenciner-Montgomery [2000], Montgomery [1998]) 
$$K = \dot R^2 + {{R^2} \over 4} 
( \dot \phi ^2 + cos^2
(\phi)
 \dot \theta ^2) + |J|^2/R^2 + \| P \|^2/M. $$
This decomposition of $K$ sometimes goes under the name
of Saari's    decomposition.
The first term $\dot R^2$ represents dilational
kinetic energy.  The last two terms   represent 
the kinetic energy of rotation and of translation.
$J$ is the total angular momentum.  $P$
is the total linear momentum.  $M$ the total mass.  The
second term  $(R^2 /4) ( \dot \phi ^2 + cos^2 (\phi)\dot \theta
^2)$ of
$K$ represents  deformations of 
 the similarity class
of the triangle.  Let us write
$$K_{shape}  =   ( \dot \phi ^2 + cos^2 (\phi)\dot \theta
^2)$$ so that this second, ``pure shape''  part of $K$ is
$(R^2 /4 ) K_{shape}$.  $K_{shape}$
  corresponds to twice the kinetic energy
of a free particle on a  unit sphere.
That sphere is  the {\it shape sphere},
the sphere whose points represent
oriented similarity classes of triangles. 

The negative of the potential $U$  can be expressed as
$$U = \tilde U( \phi, \theta)/R$$
where $\tilde U$ is a function on the sphere.

To obtain the three-body 
equations in the case of angular momentum zero,
we  set $P$ and $J$ to zero, and compute
the resulting  Euler-Lagrange equations.

\vskip .3cm

{\bf 8. Proof of theorem 2 in the case of equal masses.}

We proceed with the proof of
theorem 2 in the equal mass case.
What makes this case special is that
it is the only mass distribution for which  
 the Lagrange points coincide with the 
North and South poles of  the shape sphere.  Then the height  
$$z = sin(\phi)$$
above the equator 
is the variable of theorem 2,
where $R, \phi, \theta$ are the spherical
shape coordinates of the previous paragraph.
 The Lagrangian for the zero-angular
momentum motion is 
$$L_C = (1/2)  \dot R^2 + {{R^2} \over 4} ( \dot \phi ^2 + cos^2
(\phi)
 \dot \theta ^2) +{1 \over R} \tilde U (\phi, \theta)$$
 The  \EL equations for $\phi$ are
 $\ddt{} (\dd{L}{\dot \phi} ) = \dd{L}{\phi} ,$
or
$$\eqalign{
\ddt{}({R^2 \over 4} \dot \phi) &=
 -{R^2 \over 4}
sin(\phi) cos(\phi) \dot \theta ^2 + {1 \over R}\dd{\tilde
U}{\phi}
\cr
&= -z\{{R^2 \over 4} cos(\phi) \dot \theta ^2 -
{1 \over {R sin(\phi) }}{{\partial \tilde
U}\over {\partial \phi}}\}.}
$$  And $\dot z
= cos(\phi)
\dot
\phi$ so that 
$\ddt{}({R^2 \over 4} \dot z)
= cos(\phi) \ddt{}({R^2 \over 4} \dot \phi) + ({R^2
\over 4} \dot \phi)\ddt{cos(\phi)} 
= cos(\phi) \ddt{}({R^2 \over 4} \dot \phi) - 
{R^2 \over 4} sin(\phi ) \dot \phi^2.$
Combining this  equation with
the previous one 
and looking back at the expression
for $K_{shape}$ 
yields:
$$\ddt{}(R^2  \dot z) = -q z , 
$$ 
where
$$q  =  R^2 K_{shape} - 4 {{cos(\phi)} \over
sin(\phi)}\dd{U}{\phi}.
$$
We must  show that
$q \ge 0$, with $q = 0$ if and only
if we are at the Lagrange shape $z = \pm 1$,
with the velocity 
$(\dot R, \dot \phi, \dot \theta)$ 
satisfying  $\dot \phi =
\dot
\theta = 0$. 
   Clearly
$$K_{shape} \ge 0$$
with equality if and only if 
$\dot \phi = \dot \theta = 0$.
It remains to show that  
$$- {{cos(\phi)}\over{sin(\phi)}}\dd{U}{\phi} \ge 0
\hskip 1cm  $$
with equality if and only if $z = \pm 1$.
 We postpone
the proof of the last inequality since we will need  it for any mass distribution, 
and our proof will be 
independent of mass distribution. See (INEQ2) 
and its proof below.

\vskip .3cm

{\bf 9. Conformal geometry of the shape sphere; 
height variables.}

The variable $z$ of theorem 2 is a function on
the shape sphere.  
 The two key properties of this
variable $z$  which we  used in the proof of   theorem 1  
  are  that its zero locus 
is the equator  of collinear  configurations,   and
  that   its  critical
points   are  the Lagrange points.    
 In the equal mass both  properties
are satisfied by the height function above the equator,
$z = \sin(\phi)$, where $\phi$ is
the signed distance of a point on the sphere  from the
equator.  The North and South poles (the points a maximal distance from the equator) 
of the
shape sphere coincide with the   Lagrange
points if and only if all the masses are equal.  Consequently,
the height function above the equator fails to satisfy
the second key property in the case of unequal masses,
and we  are forced to make another choice of
the variable  $z$. 

In the case of general masses,
we take $z$ to be the  
height function   {\it as it would be defined if all the masses
were equal}.  This variable  satisfies
the two key properties, but  complicates 
the kinetic energy of the  Lagrangian.
We must  understand this complication.
The crux of the matter
is that this choice of $z$
is tantamount to applying a conformal
transformation to the  
shape sphere which takes the Lagrange points
to the North and South poles, while mapping the equator to
itself.   This conformal transformation arises via a
canonical conformal transformation from the $m$-sphere to
the equal mass distribution sphere. 

The shape space is defined to be the space
of oriented congruence classes of triangles,
while the shape sphere is the space of oriented similarity
classes of triangles. 
In other words, shape space  is the quotient of the 
three-body configuration space $(\R^2)^3$ by the group of orientation
preserving isometries,
while the shape sphere is the
quotient of $(\R^2)^3 \setminus \{
\hbox{triple collisions} \}$ by the group of orientation preserving 
similarity transformations.    As topological
spaces, neither space depends on 
the choice of masses.  The shape space
is homeomorphic to Euclidean
three space, while the shape sphere
is   homeomorphic to a two-sphere.

The triple collisions get mapped to a distinguished point of shape space,
called the triple collision point, or origin.  The action of dilation fixes this point,
while changing all  other points of the shape space.  The
shape sphere can be canonically viewed as 
the shape space minus this triple collision divided by
the action of dilations.

A choice  $m = (m_1, m_2, m_3)$ of masses defines a kinetic
energy metric on the three-body configuration
space.  This in turn induces a metric on  the shape space,
since the shape space is the quotient of the configuration
space by a group of isometries.  
The shape sphere can be realized   as the set of
 all points in  shape space
a distance $1$ from triple collision, and from here
the shape sphere inherits a metric as well. 
We denote this metric by  $d^2 s_m$. 
The shape sphere with this metric is isometric
to the standard round metric on a sphere of radius $1/2$
in Euclidean space.   We then
have that the metric on shape space is given
by 
$$dR^2 +  (1/2)^2 R^2 d^2 s_m.$$.
This  expression accounts for the kinetic energy
of the previous section. 
     
The shape sphere has a 
conformal structure which is   independent of the
kinetic energy,, i.e. is independent of the   mass distribution. 
This conformal structure is implicit in the
work of Albouy-Chenciner [1998].
We will need  the explicit   conformal factor $\lambda$
 relating   two kinetic energy metrics
on the sphere.

\proclaim Proposition.
The shape metrics
$d^2 s_m$ and $d^2 s_{m^{\prime}}$ 
for two different mass distributions  $m$ and
 $m^{\prime}$ are conformally related  according
to the formula    
$$ { {m_1 + m_2 + m_3} \over {m_1 m_2 m_3}} I_m ^2
d^2 s_m = { {m_1 ^{\prime}  + m_2 ^{\prime} + m_3 ^{\prime}} \over {m_1 ^{\prime} m_2^{\prime} m_3 ^{\prime}}} 
I_{m^{\prime}}
^2d^2 s_{m^{\prime}}$$

We will take for coordinates
on the shape sphere  standard spherical coordinates $\phi, \theta$
    for the equal mass distribution
$m ^{\prime} = (1,1,1)$ metric.
Thus   $d^2 s _{m ^{\prime}} = 
 d \phi^2 + \cos(\phi)^2 d \theta ^2$.
When we write the metric for $d^2 s_m$
in these coordinates we get 
$d^2 s_m = \lambda(\phi, \theta) (d \phi^2 +
\cos(\phi)^2 d
\theta ^2)$ with $\lambda = c(m) I_m ^2/ c(m ^{\prime}) I_{m^{\prime}}^2 $ as in the theorem,
where $c(m)$ is the total mass divided by the product of the masses. 
Recalling that the  metric defined
by the mass distribution $m$  on the three-dimensional shape space is  
$dR^2 + (R^2 /4)d^2 s _m$
where $R^2 = I_m$, 
we see that   the
kinetic energy  on shape space, 
which is obtained by setting  the total linear
and angular momentum to be zero ($P = J = 0$ in
the expression for $K$ of the previous section)  is  
$$K = \dot R^2 +   {{R^2} \over 4} K_{shape}$$ 
with 
$$K_{shape} = 
(\lambda (\phi,
\theta)( \dot \phi ^2 + cos^2 (\phi 
 \dot \theta ^2) \hskip .5cm (2).$$

The proposition implies  that the shape sphere  has a fixed
conformal structure, independent of choice of masses.  The
group of orientation-preserving conformal automorphisms of
the sphere is  the same as the group of 
orientation-preserving, circle-preserving
transformations. Thus it makes sense to 
speak of circles on  the shape sphere   without  specifying  any
  mass distribution.  

\proclaim Lemma [on circles]. Write $s_i = r_{jk} ^2$
where $ijk$ is a permutation of $123$
for  the  squared side lengths  of a triangle.  
And write $\Delta = {  1\over 2} {\bf n} \cdot ({\bf x_2 -
x_1})
\times ({\bf x_3 -x_1}) $ for its signed area.  Then the
linear equation $A s_1 + B s_2 + C s_3 + D \Delta = 0$
with $A, B, C, D$ real constants, describes
a circle in the shape sphere, provided the set of
triangles satisfying the inequality  is nonempty. 
Conversely, every circle in the shape sphere is described by
such an equation. 

The proofs of proposition and the lemma
are postponed to after the proof of theorem 2.

\vskip .3cm
{\bf 10.  Proof of theorem 2, unequal mass case.}

The proof begins  by  computing
the Euler-Lagrange equations in
our special coordinates. The computation is as for the
equal mass case, the main difference being
the occurence of $\lambda$ in the Lagrangian. 
 We compute the
Euler Lagrange equations for $\phi$,
and then for $z = \sin(\phi)$.  
We have  Lagrangian $L = (1/2)K + U$
where $K$ is given by equation (2) above.
The Euler-Lagrange equation for $\phi$ is then
$$\ddt{}({R^2 \over 4} \lambda \dot \phi )
= {1 \over 2} {R^2 \over 4} {{\partial \lambda} \over
{\partial \phi}} 
(\dot \phi ^2 + \cos(\phi)^2  \dot \theta ^2)
- {R^2 \over 4} \lambda \cos(\phi) \sin(\phi) \dot \theta^2
+ {{\partial U} \over
{\partial \phi}}$$
Using this equation and 
$z = sin(\phi)$, so that $\dot z = \cos(\phi) \dot \phi$
as in the equal mass computation,
and expanding out $\ddt{}(R^2  \lambda \dot z )$
yields
$$\ddt{}(R^2 \lambda  \dot z) = -q z , 
$$
where
$$q  = (1 -  {1 \over 2}{{\cos(\phi)} \over  {\sin(\phi)}}
{1 \over \lambda} {{\partial \lambda} \over {\partial \phi}}
) R^2 K_{shape}
 - 4{{\cos(\phi)} \over  {\sin(\phi)}} {{\partial U} \over
{\partial \phi}} \hskip .5cm (3). 
$$  
Now 
$$K_{shape} \ge 0$$
with equality if and only if all the kinetic
energy is in the dilational ($\dot R$) motion.
To conlcude the proofs then,
we require that
$$(1 - {{\cos(\phi)} \over  {\sin(\phi)}} {1 \over 2} {1 \over \lambda}
{{\partial \lambda} \over {\partial \phi}} ) > 0
\hskip .5cm (INEQ1)
$$
and 
$$ - {{\cos(\phi)} \over  {\sin(\phi)}} {{\partial U} \over
{\partial \phi}} > 0  \hskip .5cm (INEQ2)$$
 for $0 < \phi <
\pi/2$, and for the mass distribution as given. 

Note that both $U$ and $\lambda$ are even functions
of $\phi$ by reflectional symmetry. The
derivative of any function $f(\phi, \ldots)$ 
which is an even function of
$\phi$ must be  zero at $\phi = 0$,
and consequently 
${1 \over \phi}{{\partial f} \over {\partial
\phi}}$ is smooth through $\phi =0$.
It follows that both   ${{\cos(\phi)} \over 
{\sin(\phi)}}{{\partial \lambda} \over {\partial \phi}}$ and
${{\cos(\phi)} \over  {\sin(\phi)}}{{\partial U} \over
{\partial \phi}}$ are smooth functions through the equator. 

\vskip .2cm
{\bf Proof of Inequality 2.}
The inequality (INEQ2) is valid for all mass distributions.
Since $cos(\phi)/sin(\phi)$ is an odd function,
positive for $0 < \phi < \pi/2$,
and since ${{\partial U} \over
{\partial \phi}}$ is also odd, it suffices
to show that $- {{\partial U} \over
{\partial \phi}}$ is positive in the range 
$0 < \phi < \pi/2$.

The proof of the positivity of
$-{{\partial U} \over
{\partial \phi}}$
is elegant but tricky.
Introduce as coordinates in shape space
$$s_k = r_{ij} ^2$$
for $ijk$ any permutation of $123$.  Then
$$U =   
m_1 m_2 /s_3 ^{1/2}  + 
m_3 m_1 /s_2 ^{1/2} + m_2 m_3 /s_1 ^{1/2} $$
while
$$I = (m_1 m_2 s_3 + m_3 m_1 s_2 + m_2 m_3 s_1) /M \hskip
.5cm (L1).$$
with $M = m_1 + m_2 + m_3$.  
To differentiate with respect to $\phi$
we fix $I$ and $\theta$, thus defining meridianal
circles passing through the Lagrange point,
and then differentiate along these meridianal curves.
The crux of the inequality is 
to observe that each of these meridianal curves 
is defined by a linear constraint
$$A s_1 + B s_2 + C s_3 = 0
\hskip
.5cm (L2)$$
 when written in
terms of the $s_k$.
Here   $A, B,C$ are any real constants, not all zero,
but summing to zero.
To see the validity of this representation
of the meridianal curves, use the lemma
of the previous section. It says that any circle  in the
shape sphere can be expressed in the
form $A s_1 + B s_2 + C s_3 + D \Delta = 0$.
Now the   meridianal circles  pass through
 the two Lagrange points $L_+$ and $L_-$,and
any circle passing through these two points is
a meridianal circle.   The Lagrange points are characterized
by
$s_1 = s_2 = s_3$,  while their signed areasare
are negatives of each other:
$\Delta (L_+ ) = - \Delta (L_-)$.
Writing $s_i =s$ and $\Delta = \Delta( L_+)$
we see that the   
the coefficients defining the circles satisfy 
$(A + B + C)s + D \Delta = 0$
and $(A + B + C)s -  D \Delta = 0$.
Neither $s$ nor $\Delta$ are zero. 
Subtracting  the two equations yields $D = 0$.
Adding them yields $A + B + C = 0$.

Since $1/s^{1/2}$ is convex for
$s >0$,   $U$ is
a  strictly convex function in the positive coordinate
orthant $s_k > 0$.  The constraints 
(L1) and (L2) are linear,
so upon restriction, $U$ is again a strictly convex
function.  Consequently,
with the constaints imposed, $U$  has at most one global
minimum.
But (either of) the  Lagrange point 
$L$ (i.e  $L_{+}$ or $L_-$) is the   global minimum of $U$
when we  impose only constraint (L1). 
(Note that $\Delta$ does not occur in the 
constraints or in the expression for $U$.
In essence we are also allowing reflections when
we ignore $\Delta$ and use only the $s_i$ as
coordinates on the shape space.)
  All the lines
defined by (L2) pass through $L$. 
Consequently, 
$U$ restricted to the
line (meridian) defined by both constraints (L1) and (L2)  
has a unique minimum at $L$
and    is strictly
increasing as we move away from it.  The variable
$\phi$
 monotonically decreases as we move away
from $L$ toward the equator.  This proves that
$${{\partial U} \over {\partial \phi}} <0$$
for all $\phi$ with $0 < \phi < \pi/ 2$.

\vskip .2cm
{\bf Proof of Inequality 1.}

We can rewrite the desired inequality as
$$1 + {\cos(\phi) \over \sin(\phi)} \dd{}{\phi} log{\hat I}$$
where 
$$\hat I := I/ I_1$$
and where I have used the fact 
$\lambda = C I_1 ^2/ I^2$ so that $-{1 \over 2}{1 \over
\lambda} 
\dd{}{\phi} \lambda = + \dd{}{\phi} log{\hat I}$.

To compute this logarithmic  derivative of $\hat I$,
define variables 
$$\hat s_i := s_i/I_1 = r_{jk}^2/I_1, $$
so that
$$\hat I   = {1 \over M} m_1 m_2 \hat s_3 + m_3 m_1 \hat s_2
+ m_2 m_3 \hat s_1,$$
 where $M = m_1 + m_2 + m_3$.
We need to be able to differentiate
$\hat s_i$ with respect to $\phi$.
This is easy once we  have the representation: 
$$\hat s_i = 1 - cos(\phi) \gamma_i (\theta) 
\hskip .5cm (3)$$
which we now explain, following Chenciner-Montgomery
[2000], pp. 890-891, or the end of the appendix here.

We can represent a point in shape space
as a 3-vector ${\bf w}$ in Euclidean 3-space which
we express in spherical coordinates as 
$${\bf w} = I_1( cos (\phi) \cos (\theta) , \cos(\phi)
\sin(\theta), \sin(\phi)).$$ 
Then $I_1 = \| {\bf w} \|$, while  
$\Delta = I_1 \sin(\phi) =$ a signed area,
and 
$$s_k : = r_{ij}^2 = |{\bf w}| - {\bf w} \cdot {\bf b}_k,$$
where $ijk$ is a  permutation of $123$, and where the 
 ${\bf b}_k$ are three unit vectors on the equator $\phi = 0$
which represent the binary collision rays.  
These three unit vectors ${\bf b}_k$ are arranged at the
vertices of an equilateral triangle circumscribed in the
unit circle.  Write
${\bf u} = (\cos(\theta), \sin(\theta), 0)$
 and 
$$\gamma_k (\theta):= {\bf u}\cdot {\bf b}_k.$$ 
Then we have  
that $s_k = I_1 - I_1 \cos(\phi) \gamma_k (\theta)$
and the equation (3) for the $\hat s_i$ follows immediately.

Writing
$$p_k = m_i m_j/M >0 $$  
we have 
$$\hat I = \Sigma p_i \hat s_i.$$
Using the expression (3) for
$\hat s_i$ we  compute:
$$\dd{}{\phi} log{\hat I} = \Sigma p_k \sin(\phi) \gamma_k/
\Sigma p_k \hat s_k.$$
It follows that
$${\cos(\phi) \over \sin(\phi)} \dd{}{\phi} log{\hat I}
= \Sigma p_k \cos(\phi) \gamma_k/
\Sigma p_k (1 - \cos(\phi) \gamma_k, $$
and 
$$1 +  {\cos(\phi) \over \sin(\phi)} \dd{}{\phi} log{\hat I}
= \Sigma p_k / \Sigma p_k ( 1 - \cos(\phi) \gamma_k).$$
Now  use the fact that 
$|\cos(\phi) \gamma_k | \le | \gamma_k | \le 1$ and that at
least one of the $|\gamma_k|$ is less than $1$ to conclude
that the previous expression  is finite and positive.

QED

\vskip .5cm

{\bf 11. Proofs of the proposition and
the lemma ;Conformal Geometry.}

We will give two different proofs of
proposition, and one proof of the lemma.

\vskip .2cm

{\bf 11.1.  Proof of the proposition via  Jacobi coordinates.}

Write $E = \R^2 \times \R^2 \times \R^2$
for the configuration space of the
three-body problem.
The $ith$  Euclidean plane
 factor represents the  positions of the
$i$th  body. Write points of $E$ as
$x= (x_1, x_2, x_3) \in E$
with $x_i \in \R^2$.
Identify $\R^2$ with the complex numbers $\C$
in the standard way so that $E = \C^3$.
The Jacobi map ${\cal J}_m$
associated to the mass distribution $m = (m_1, m_2, m_3)$
is the linear map 
$${\cal J}_m: E \to \C^2$$ given by 
$${\cal J}_m (x_1, x_2, x_3) = (z_1, z_2)$$ 
where 
$$z_1 = \sqrt{\mu_1} (x_2 - x_1),$$
$$z_2 = \sqrt{\mu_2} (x_3 - ((m_1 x_1 + m_2 x_2)/(m_1 + m_2))$$
and
$${1 \over {\mu_1}} = {1 \over m_1} + {1 \over  m_2},$$
$${1 \over {\mu_2}} = {1 \over m_3} + {1 \over  {m_1 +
m_2}}.$$ 
Physically $z_1$ is the normalized edge vector
joining 1 to 2, and $z_2$ 
is obtained by normalizing the vector 
which joins the center of
mass of this edge to the remaining vertex.  

The Jacobi map
 is invariant under translations: 
${\cal J}_m ((x_1 + v, x_2 + v, x_3 + v) = {\cal J}_m (x_1,
x_2, x_3)$.
It  diagonalizes the kinetic energy 
$$\eqalign{
K & :=m_1 \| \dot x_1 \|^2 + m_2 \| \dot x_2 \|^2 +
m_3
\|
\dot x_3
\|^2
\cr
& = \| \dot z_1 \|^2 + \|\dot z_2 \|^2  
}$$ 
provided the total linear momentum is zero:
$m_1  \dot x_1 + m_2  \dot x _2 + m_3 \dot x_3  = 0$.
Similarly, it diagonalizes the moment of inertia
tensor:
$$\eqalign{I & := 
m_1 \|  x_1 \|^2 + m_2 \|   x_2 \|^2 +
m_3 \|  x_3\|^2
\cr
& = \|   z_1 \|^2 + \|  z_2 \|^2  
},$$
provided the center of mass is at the origin
$m_1   x_1 + m_2   x _2 +   m_3 x_3   = 0$. 

The action of the group of orientation
preserving similarities on triangles
$x$ becomes, under the Jacobi map,
the action of complex scalar multiplication:
$(z_1, z_2) \mapsto (\lambda z_1, \lambda z_2)$,
$\lambda \in \C$, $\lambda \ne 0$.  
Thus the shape sphere is identified with the
complex projective line $\C \bP ^1$, the 
space whose points are complex lines in $\C^2$. 
  The quotient map  
$$\pi: \C^2 \setminus \{(0,0) \} \to \C P^1 = S^2$$
sends a nonzero  complex
vector $(z_1, z_2)$ to the 
  complex line $\pi(z_1, z_2)= [z_1, z_2]$ which 
it spans.   The map $\pi \circ {\cal J}_m: E \setminus 
\{ \hbox{triple
collisions} 
\} \to \C \bP ^1 = S^2$
 sends a triangle $x \in E$
to its ``shape'' meaning oriented
similarity class. Note that  
we must delete the triple collisions 
$x_1 = x_2 =x_3$ because they form the  
kernel of the Jacobi map.

If we now repeat the procedure with a different 
mass distribution $m^{\prime} = (m_1 ^{\prime}, 
m_2 ^{\prime}, m_3 ^{\prime})$ we obtain
different Jacobi coordinates  $w_1, w_2$,
which diagonalize the new  moment of inertia
$I_{m^{\prime}}$.

\vskip .2cm

We  abstract the situation described above. 
Consider  a complex
two-dimensional vector space $\C^2$ with its
standard complex structure.   This
vector space represents the space of
Jacobi coordinates. Write $\C \bP ^1$ for the
corresponding complex projective line. It is the  quotient of $\C^2
\setminus \{ 0\}$ by the action of complex scalar
multiplication.

A Hermitian inner product on $\C^2$ induces 
a metric on $\C \bP ^1$ as follows.  
Write $I (z) = \langle z ,z \rangle$
for   square norm for this Hermitian inner product.
Setting $I =1$ defines a  three-sphere $S^3 _I$
with induced \Ri metric coming from the real part of the 
Hermitian innerproduct. .  
The subgroup $S^1 \subset  \C^*$ preserves $I$,
and  the inner product, and hence acts
on $S^3 _I$ by isometries.  
Consequently the quotient $S^3 _I /S^1$
inherits a \Ri metric by declaring the
submersion $S^3 _I \to S^3 _I /S^1$
to be a \Ri submersion.   The quotient
space $S^3 _I/ S^1$ is canonically identified
with $\C \bP ^1$ by sending the $S^1$-orbit of
a point $z \in S^3 _I$ to the corresponding
$\C^{\ast}$ orbit.  In this way, we obtain
a \Ri metric $d^2 s _I$ on $\C \bP ^1$.
If $(z_1, z_2)$ are  Hermitian orthonormal coordinates
so that $I = |z_1 |^2 + |z_2|^2$,
and if $z = z_1/z_2$ are the corresponding
affine coordinate on $\C \bP ^1$,
then 
$$ds _I = |dz|/(1 + |z |^2) \hskip .5cm \hbox{ for } \hskip .5cm  I = |z_1|^2 + |z_2|^2.$$
Consider another Hermitian inner product,
with corresponding square norm $I^{\prime}$.
We then have   another metric
$$ds_{I^{\prime}} = |dw|/(1 + |w|^2)
 \hskip .5cm \hbox{ for } \hskip .5cm  I^{\prime} = |w_1|^2 + |w_2|^2$$
on the same projective space, 
but now with   affine coordinate
$w = w_1/w_2$. The proposition becomes
a special case of 
\proclaim Theorem 3.  Let $I$ and $I^{\prime}$
be the square norms for two
different Hermitian structures on
the same complex two-dimensional vector space.
Let $\C \bP
^1$ be the projectivization of this vector space,
and let $ds_{I}$ and $ds_{I^{\prime}}$ 
be the two metrics on this projective space induced
by our two Hermitian inner products.
Let $L: V \to V$ be a linear operator
intertwining the two norms:  $I(Lz) = I^{\prime} (z)$.
Then the  two metrics  are related by 
$$ds_{I^{\prime}} = |det(L)| (I/I^{\prime}) ds_{I}.$$

{\bf Proof of theorem 3.} From basic linear algebra,
the  complex linear intertwining map $L$ of the theorem
always exists.  It is found by choosing orthonormal
coordinates $(z_1, z_2)$  for $I$, expressing the inner product
for $I^{\prime}$ as a matrix
in these coordinates, and then diagonalizing this matrix. 
If
$$L = \pmatrix{ a & b \cr c & c }.
$$
then
$$w_1 = a z_1 + b z_2$$
$$w_2 = c z_1 + d z_2$$
are orthonormal coordinates for the Hermitian inner product
with square norm $I^{\prime}$.
 The corresponding affine coordinates $z = z_1/z_2$
and $w = w_1/w_2$  are then 
related by the   linear fractional
transformation
$$w = (az + b) /(cz + d).$$
We compute
$$dw = (ad - bc) dz/(cz + d)^2.$$
(We ask our  gentle reader to please  bear
with us and not be confused by the two meanings
of the letter ``$d$'' here.)  Setting 
$D = |ad -bc| = |det(L)|$, we have 
$$\eqalign{
{|dw| \over {1 + |w|^2}} 
& = {{1 + |z|^2} \over {1 + |w|^2}}
{D |dz| \over {|cz + d|^2}} {1 \over {1 + |z|^2}} \cr
& = {{1 + |z|^2} \over {|cz + d|^2 + |a z + b|^2}}
{D |dz| \over {1 + |z|^2}} \cr
& = {{|z_2|^2 + |z_1|^2} \over {|cz_1 + dz_2|^2 + |a z_1 +
bz_2|^2}} {D |dz| \over {1 + |z|^2}} \cr 
& =  {I \over I^{\prime}} D { |dz| \over {1 + |z|^2}}
}.
$$
In the third line we multiplied both the numerator and
denominator of the first fraction by $|z_2|^2$. 
 QED

\vskip .2cm

{\bf Completion of the Proof of the proposition.}
Theorem 4 tells us that
$d ^2 s_{m^{\prime}} =  C (I_m ^2) / (I_{m ^{\prime}}^2) d^2 s_m$
and that the constant $C$ is given by 
$C = |det(L)|^2$ where 
 $L$ is an intertwining operator taking $I_m$ to 
$I_{m^{\prime}}$. To complete the proof of the proposition
we  solve for   $L$
so as to obtain the   correct constant $C$.  

Fix the triangle 
$x = (x_1, x_2, x_3) \in E = \R^2 \times \R^2 \times \R^2$,
 the configuration space of the
three-body problem. Then
it has two images $z$ and $w$ in $\C^2$ according to
the Jacobi maps for the two different   mass distributions
$m$, and $m^{\prime}$. 
Write 
$z=  {\cal J}_{  m} (x)$
and
$w= {\cal J}_{m^{\prime}} (x)$.

We look for a linear map
$L: \C^2 \to \C^2$ such that
$w = Lz$.  Make the upper triangular anzatz 
$L(z_1, z_2) = ( \alpha  z_1, \beta z_1 + \gamma z_2)$.
Using the above expression
for the Jacobi map, the ansatz leads to 
the  two linear equations
$\alpha z_1 = w_1$ and $\beta z_1 + \gamma z_2 = w_2$,
or
$$ \alpha   \sqrt{\mu_1} (x_2 - x_1) = \sqrt{\mu_1 ^{\prime}} (x_2 -
x_1),$$ and 
$$\beta \sqrt{\mu_1} (x_2 - x_1) + 
\gamma  \sqrt{\mu_2} (x_3 - (m_1
x_1 + m_2 x_2)/(m_1 + m_2) )
=  \sqrt{\mu_2 ^{\prime}} (x_3 - (m_1 ^{\prime}
x_1 + m_2 ^{\prime} x_2)/(m_1 ^{\prime} + m_2 ^{\prime})).$$
The first equation has 
$\alpha = \sqrt{{\mu_1 ^{\prime}}/{ \mu_1} }$
for a solution. 
Expanding out the second  equation in $x_1, x_2,x_3$
and equating coefficients yields a system of 
three homogeneous equations, in the two unknowns $\beta$ and
$\gamma$.  
 The $x_3$
equation has
$\gamma = \sqrt{\mu_2 ^{\prime} / \mu_2}$ as a solution.
Using this $\gamma$, the $x_1$
equation has 
 $\beta = -\sqrt{\mu_2 ^{\prime} /\mu_1}
( m_1/ (m_1 + m_2) - m_1 ^{\prime}/ (m_1 ^{\prime} + m_2
^{\prime}))$ for a solution,  while the $x_2$  equation 
 $\beta = \sqrt{\mu_2 ^{\prime} /\mu_1}
(m_2/ (m_1 + m_2) - m_2 ^{\prime}/ (m_1 ^{\prime} + m_2
^{\prime}))$ has for   solution. These two $\beta$s are  checked to be
equal, and so we get our invertible  linear operator 
$$L = \pmatrix{ \alpha & 0 \cr \beta & \gamma }.$$
We have $det(L) = \alpha \gamma = \sqrt{ \mu_1 ^{\prime} \mu_2 ^{\prime}/\mu_1 \mu_2}$.
Plugging in the formulae for the $\mu$ in terms of the masses
leads to $\mu_1 \mu_2 = m_1 m_2 m_3/ (m_1 + m_2 + m_3): = c(m)$.
Consequently $det(L)  = \sqrt{c(m^{\prime})/c(m) }$. 
Finally, plugging in $m_{\prime} =  (1,1,1)$ yields the formula
of the proposition.

\vskip .3cm

{\bf 11.2.  Invariant theory.} 

In order to obtain another proof of Theorem 3, we  search  for
a metric-independent  
geometric interpretation of expression 
$I_m ^2 d^2 s_m$.   This alternative point of view
will also yield a simple proof of the lemma on circles. 

Consider the vector space $V$ of planar triangles
 modulo translation, i.e.   
$(\R^2)^3$ modulo translations. 
$V$ is a complex two-dimensional
vector space, which is to say a
 real vector space  endowed
with an almost complex structure $J$, but with no
 canonical  inner product.  The inner product
must await  the introduction of masses.   
$J$ rotates triangles by ninety
degrees counterclockwise.   The circle group $S^1$ acting on triangles
by rotation   consists of the transformations  $exp(\theta J)$,
 $\theta$ real.   

Consider the real vector space ${\cal P}$ 
 of real  quadratic  $S^1$-invariant polynomials on
$V$ which are invariant
under the  action of the circle group.  ${\cal P}$ is also a
four-dimensional real vector space.   
One choice of basis for ${\cal P}$
consists of the
the squared side lengths $s_k =r_{ij}^2$
and the signed area $\Delta$.  
Another choice of basis is obtained by
choosing complex linear coordinates, for example
Jacobi coordinates,  
$z_1, z_2$   
for $V$.  Then  
$|z_1|^2, |z_2|^2$ and the real and imaginary parts of $z_1
\bar z_2$ form a   basis for  ${\cal P}$.
If  
$\langle z, w \rangle = z_1 \bar w_1 + z_2 \bar w_2$
denotes the standard Hermitian form relative to these coordinates,
then  we can identify ${\cal P}$ with the space 
${\cal H}$ of two-by-two Hermitian
matrices.  For any invariant $I$ can be expressed uniquely in the form  
$$I(z) = \langle z , H z \rangle$$
for some unique Hermitian matrix   $H \in {\cal H}$. 

Every $S^1$-invariant function $f$ is expressible as
a function in the quadratic invariants.
It follows that if we know the values of a point $v \in V$
on a basis for ${\cal P}$, then we know the $S^1$-orbit of $v$.
Let ${\cal P}^*$ be the vector space dual
to ${\cal P}$.  For  $v \in V$, 
define a linear functional $ev(v)$,
the {\it evaluation map},  on ${\cal P}$
by:
$$ev(v)(Q) = Q(v).$$
This evaluation map is a canonical map 
$$ev: V \to {\cal P}^*.$$
and according to what we have just said, its
image is a realization of the quotient space
$V/S^1$, i.e. of ``shape space''.
\proclaim Lemma.  The image $ev(V)$ of the evaluation map
is isomorphic to the quotient space $V/S^1$.
This image is the positive half of a quadratic cone
in the vector space ${\cal P}^*$, the cone being defined
by the vanishing of a real quadratic form of signature $(3,1)$.
Consequently,  ${\cal P}^*$ and ${\cal P}$  are endowed with
canonical Minkowski inner products,  denonted $\beta(v, w)$, unique up to scale.

A choice of basis for ${\cal P}$ is a system of linear coordinates
on ${\cal P}^*$.  The  cone of the lemma can be described as
a  quadratic relation between the elements of the  basis.
If we choose for basis  the   squared side-lengths   $s_k = r_{ij}^2$, 
together with  the signed area $\Delta$ of the triangle, then the cone 
results from  Heron's  relation 
$$16 \Delta^2 = ( r_{12} + r_{23} + r_{31} ) 
( r_{12} + r_{23} - r_{31} ) 
(r_{23} + r_{31} -  r_{12} )
 (r_{31} + r_{12} - r_{23} ).$$
Expand the right hand side to obtain
$$16 \Delta ^2 = 2 s_1 s_2 + 2 s_3 s_1 + 2 s_2 s_3 - (s_1 ^2 + s_2 ^2 + s_3^2)
\hskip .5cm ; \hskip .5cm s_i \ge 0$$
which describes the positive half of the cone of the lemma.
 If instead we use  the basis  $|z_1|^2, |z_2|^2$, $Re(z_1 \bar z_2), Im(z_1 \bar z_2)$
then the  cone results from   the relation 
$(|z_1|^2 |z_2|^2 = |z_1 \bar z_2|^2$.
Alternatively, take the basis
$w_0 =  {1 \over 2}(|z_1|^2 + |z_2|^2)$, 
$w_1 = {1 \over 2}(|z_1|^2 - |z_2|^2)$,
$w_2 = Re(z_1 \bar z_2)$,
$w_3 = Im (z_1 \bar z_2)$.  Then the positive  cone is given by
$w_0 ^2 = w_1 ^2 + w_2 ^2 + w_3 ^2$, $w_0 \ge 0$, 
a relation which holds among the functions
at all points  of $V$.
(This relation is  familiar from 
the Hopf map.) If we use the  coordinates $z_1, z_2$ to view ${\cal P}$
as ${\cal H}$, then we can also identify
${\cal P}^*$ with ${\cal H}$ using the trace pairing
to identify ${\cal H}$ with
${\cal H}^*$.  In these coordinates:
$$ev(z_1, z_2)_{ij} = H_{ij} (z): = z_i \bar z_j.$$
and the cone is defined by the relation
$$det(H) = 0 \hskip .5cm ; \hskip .5cm tr(H) \ge 0.$$

The  group 
$GL(V ; J)$ of  linear transformations of $V$
which commute with $J$ acts linearly   on the invariants by pull-back,
and hence acts linearly on  ${\cal P}^*$.  By construction,
this action  
  preserves the quadratic coneand so 
is an action by means of  
the linear conformal Lorentz group  $CSO_+(\beta) \cong CSO(3,1)_+$.
Here the subscript $+$ denotes the time orientation preserving part of 
the full Minkowski isometry group, and
the $S$ denotes the   orientation
preserving part.   If we fix a complex volume
element  in $V$, and hence restrict $GL(V ; J)$ to
$SL(V)$, the action just defined is  the well-known $2:1$ homomorphism
$SL(2, \C) \to SO(3,1)$.

Now let us projectivize, which is to say, divide by
dilations.  These dilations correspond to scaling similarities
of our triangle.  Now the set of   rays in the light cone in Minkowski space
forms a two-sphere.  This is our shape sphere.
The action of $GL(V, J)$, which factors through 
$CSO_+(\beta)$ as we have just seen,
is an action on this   sphere
by conformal transformations. 
 Now we are ready to prove the theorem 3.

\vskip .2cm

{\bf Second Proof of Theorem 3.}

Fix a representative Minkowski structure
$\beta$ on ${\cal P}^*$, one whose
cone $C = \{p: \beta(p, p) = 0 \}$
is our quadratic cone.  The restriction 
$\beta_C$ to the cone is a degenerate metric of signature $(2,0)$.
If $(x,y,z,t)$ are standard Minkowski orthonormal
coordinates for $({\cal P}^*, \beta)$ then
$\beta = dx^2 + dy^2 + dz^2 - dt^2$
while $C$ is defined by $x^2 + y ^2 + z^2 -t^2 = 0$.
Write $r^2 = z^2 + y^2 + z^2$.  
Write $d^2 \sigma_t$ for  the restriction
of $\beta$ to the two-sphere $r = 1$  in 
the space-like hyperplane $t = 0$. 
We compute
$$\beta_C = r^2 d^2 \sigma_t = t^2 d^2 \sigma_t.$$

More generally, if $\tau$ is any time-like linear coordinate
then
$$\beta_C = {{\tau ^2} \over {(\tau, \tau)}} d^2 \sigma_{\tau}.$$
where the numerical constant $(\tau, \tau)$ is the Minkowski
length of the dual vector $\tau \in {\cal P}$.  
To see this, write $t = c \tau$ where $(\tau, \tau) = 1/c$,
thus defining a unit time-like linear coordinate
which can be completed to form a system  
$(x,y,z,t)$ of  Minkowski orthonormal
coordinates. In this formula, $d^2 \sigma_{\tau}$ is again
the restriction of $\beta$ to the unit sphere in
the space-like Euclidean hyperplane $\tau = 0$.

The square norm $I$ for any Hermitian inner-product on
$V$ is a linear time-like coordinate on ${\cal P}^*$.
Thus if $I, I^{\prime}$ are two such square norms we have:
$${{I^2} \over {(I, I )}} d^2 \sigma_{I} = \beta_C =
{{I^{\prime 2}} \over {(I^{\prime}, I^{\prime} )}} d^2 \sigma_{I^{\prime}}.$$
We are almost done. It remains to evaluate the constant
$(I^{\prime}, I^{\prime} )/(I, I)$.  If $H$ is the Hermitian matrix representing
$I$ in some system of coordinates, then
$H^{\prime} = L H L^*$ represents $I^{\prime}$ where
$L$ is the intertwining operator.  But we have
seen that a choice for the Minkowski inner product is
$(I, I) = det(H)$, and $det(H^{\prime}) = |det(L)|^2 det(H)$,
so that $(I^{\prime}, I^{\prime} )/(I, I) =  |det(L)|^2 $.
QED

\vskip .2cm {\bf Remark.}  A choice of square norm $I$
fixes a normalization of the Minkowski inner product
$\beta$ by declaring that $(I, I) = 1$.
With this normalization, the shape space metric
on the cone is ${1 \over 4} \beta_C + dI^2$.

\vskip .2cm

{\bf 10.3.  Proof of the lemma on circles.}
Circles on a sphere are obtained by
intersecting the sphere with planes.  
Think of the sphere as the 
projectivized cone in Minkowski space.
Realize this sphere as in the second  proof
of theorem 3 by  intersecting the quadratic cone
in ${\cal P}^*$ 
with the three-dimensional affine space   $\{ I =1 \}$, 
where $I$ is the square norm for  a  $J$-compatible
inner product on $V$.  
The   $s_i$ and $\Delta$
form linear coordinates on ${\cal P}^*$,
and so by restriction any three of them form
linear coordinates on the affine space $I = 1$.
The planes in this affine space are  defined
by  a linear equation in the $s_i$ and $\Delta$.

QED

\vskip .4cm
{\bf Acknowledgements.}
I would like to thank  
Alain Albouy, Alain Chenciner, 
Mark Levi, Rick Moeckel
and  Jeff Xia  for
conversations crucial to the development
of this paper.  Most of the ideas
in 10.2, the proof of theorem 3 via
invariant theory,  are due to Albouy.
This work was supported in part by 
NSF grant  (DMS 9704763).

\vskip .4cm

{\bf Bibliography.} 

A. Albouy and A. Chenciner, [1998],
Le probl\'eme des n corps et les distances
mutuelles, Inventiones,  131, 151-184.

A. Chenciner  and R. Montgomery,  [2000], 
A remarkable periodic solution 
of the three-body problem in the case of equal masses, 
{\it Annals 
of Mathematics}, {\bf 152}, 881-901. 

 R. Moeckel, [1983], Orbits 
Near Triple Collision in the Three-Body Problem, 
Indiana Univ. Math. J., v. 32,  no. 2, 221-240.


\smallskip

R. Montgomery, [1998], The $N$-body problem, the braid
group, and action-minimizing orbits,  Nonlinearity,  
{\bf 11}, 363-376.
 

\end